\documentclass{article}

\usepackage{amssymb,latexsym,amsfonts,epsfig}

\newtheorem{theorem}{Theorem}[section]
\newtheorem{lemma}[theorem]{Lemma}
\newtheorem{alg}[theorem]{Algorithm}

\newcommand{\hpg}[5]{{}_{#1}\mbox{\rm F}_{\!#2}\!
  \left(\left.{#3 \atop #4}\right| #5 \right) }

\newcommand{\singset}{\Delta}
\newcommand{\proof}{{\bf Proof. }}
\newcommand{\qed}{\hfill $\Box$\\}
\newcommand{\equal}{\!\!=\!\!}
\newcommand{\CC}{{\Bbb C}}

\newcommand{\PP}{{\Bbb P}}
\newcommand{\HH}{{\Bbb H}}

\newcommand{\RR}{{\Bbb R}}

\title{Transformations of some Gauss
hypergeometric functions}

\author{Raimundas Vid\=unas\footnote{Supported by NWO,
project number 613-06-565, and by the Flemish FWO NOG-project.}\\
\em Kyushu University}

\date{}

\begin{document}
\maketitle

\begin{abstract} This paper presents explicit algebraic
transformations of some Gauss hypergeometric functions. Specifically,
the transformations considered apply to hypergeometric solutions of
hypergeometric differential equations with the local exponent differences
$1/k,1/\ell,1/m$ such that $k,\ell,m$ are positive integers and
$1/k+1/\ell+1/m<1$. All algebraic transformations of these Gauss
hypergeometric functions are considered. We show that apart from 
classical transformations of degree 2, 3, 4, 6 there are several other
transformations of degree 6, 8, 9, 10, 12, 18, 24. Besides, we present an
algorithm to compute relevant 
Belyi functions explicitly.
\end{abstract}

\section{Introduction}

An algebraic transformation of Gauss hypergeometric functions
is an identity of the form
\begin{equation} \label{hpgtransf}
\hpg{2}{1}{\!\widetilde{A},\,\widetilde{B}\,}{\widetilde{C}}{\,x}
=\theta(x)\;\hpg{2}{1}{\!A,\,B\,}{C}{\varphi(x)},
\end{equation}
where $\varphi(x)$ is a rational function of $x$, and $\theta(x)$
is a product of some powers of rational functions. Here are two
examples of quadratic 
transformations (see \cite{bateman,specfaar}):
\begin{eqnarray} \label{hpgtr2}
\hpg{2}{1}{a,\,b}{\frac{a+b+1}{2}\,}{\,x} & = &
\hpg{2}{1}{\frac{a}{2},\;\frac{b}{2}\,}{\frac{a+b+1}{2}\,}{\,4x\,(1-x)},\\
\label{hpgtr2a} \hpg{2}{1}{a,\;\frac{a-b+1}2}{a-b+1}{\,x} & = &
\left(1-\frac{x}{2}\right)^{-a}
\hpg{2}{1}{\frac{a}{2},\,\frac{a+1}{2}\,}
{\frac{a-b}{2}+1}{\frac{x^2}{(2\!-\!x)^2}}.
\end{eqnarray}
These identities hold in some neighborhood of $x=0$ in the complex
plane, and can be continued analytically. For example, formula
(\ref{hpgtr2}) holds for $\mbox{Re}(x)<1/2$.

Recall that the Gauss hypergeometric function
$\hpg{2}{1}{\!A,\,B}{C}{z}$ 
is a solution of the {\em hypergeometric differential equation}
\begin{equation} \label{hpgde}
z\,(1-z)\,\frac{d^2y(z)}{dz^2}+
\big(C-(A\!+\!B\!+\!1)\,z\big)\,\frac{dy(z)}{dz}-A\,B\,y(z)=0.
\end{equation}
This is a Fuchsian equation on the complex projective line $\PP^1$
with 3 regular singular points $z=0,1$ and $\infty$. 
The local exponent differences at these points are (up to a sign)
$1-C$, $C-A-B$ and $A-B$ respectively.

Algebraic transformations of Gauss hypergeometric functions usually come
from those transformations of hypergeometric equation (\ref{hpgde}), which
have the form
\begin{equation} \label{algtransf}
z\longmapsto\varphi(x), \qquad y(z)\longmapsto
Y(x)=\theta(x)\,y(\varphi(x)),
\end{equation}
and such that the transformed equation for $Y(x)$ is a
hypergeometric equation in the new indeterminate $x$. 
Here $\varphi(x)$ and $\theta(x)$ have the same meaning as in formula
(\ref{hpgtransf}). Geometrically, this is a {\em pull-back transformation}
of equation (\ref{hpgde}) with respect to the finite covering
$\varphi:\PP^1\to\PP^1$ determined by the rational function $\varphi(x)$. In
\cite{kitaevdb} these transformations are called {\em RS-transformations}.
Recall that a rational function on a Riemann surface is a {\em Belyi
function} \cite{shabatgb,kreines} if it has at most 3 critical values, or
equivalently, if the corresponding covering of $\PP^1$ branches only above a
set of 3 points. The function $\varphi(x)$ in hypergeometric identities like
(\ref{hpgtransf}) is usually a Belyi function.

Algebraic transformations of Gauss hypergeometric functions and pull-back
transformations of hypergeometric equations are related as follows.
\begin{lemma} \label{transeqv}
\begin{enumerate}
\item Suppose that pull-back transformation $(\ref{algtransf})$ of equation
$(\ref{hpgde})$ is a hypergeometric equation as
well, 
and that the transformed equation has non-trivial monodromy. Then, possibly
after fractional-linear transformations on the projective lines, 
there is an identity of the form $(\ref{hpgtransf})$ between hypergeometric
solutions of the two hypergeometric equations.%
\item Suppose that identity $(\ref{hpgtransf})$ holds in some region of the
complex plane. Let $Y(x)$ denote the left-hand side of the
identity. If 
$Y'(x)/Y(x)$  is not a rational function of $x$, then the transformation
$(\ref{algtransf})$ converts the hypergeometric equation $(\ref{hpgde})$
into a hypergeometric equation for $Y(x)$.
\end{enumerate}
\end{lemma}
\proof This is Lemma 2.1 in \cite{algtgauss}. \qed

In this paper, we consider Gauss hypergeometric functions which satisfy
hypergeometric equations with local exponent differences $1/k,1/\ell,1/m$
such that $k,\ell,m$ are positive integers and $1/k+1/\ell+1/m<1$. We call
these functions {\em hyperbolic hypergeometric functions}, because they have
interesting analytic properties related to the hyperbolic geometry of the
complex plane \cite{yoshida,beukers}. The main purpose of this paper is to
describe algebraic transformations of these functions into other
hypergeometric functions. Existence of their non-classical transformations
of degree 10, 12 and 24 is shown in \cite{hodgkins2,beukers}. A
transformation of degree 8 is presented in \cite[Section 5]{kitaevdb}. We
give a complete list of possible algebraic transformations of hyperbolic
hypergeometric functions. Algebraic transformations of all Gauss
hypergeometric functions are classified in \cite{algtgauss}.

For hyperbolic hypergeometric functions, algebraic transformations always
induce pull-back transformations of their hypergeometric equations, and vice
versa. Indeed, Kovacic algorithm \cite{kovacic,marius} in differential
Galois theory implies that the monodromy group of those hypergeometric
equations is not trivial, and that they have no solutions $y(z)$ with
algebraic logarithmic derivative $y'(z)/y(z)$. Therefore Lemma
\ref{transeqv} allows no exceptions.

This paper classifies algebraic transformations 
of hyperbolic hypergeometric functions by finding all pull-back
transformations of their hypergeometric equations to other hypergeometric
equations. The main problem is to compute suitable coverings
$\varphi:\PP^1\to\PP^1$. Possible branching patterns for them are derived in
Section \ref{otherat}. We give a general algorithm for computing coverings
with prescribed branching pattern in Section \ref{fincovers}. Algebraic
transformations of hyperbolic hypergeometric functions are listed in Section
\ref{hyperids}.

\section{Possible branching patterns}
\label{otherat}

A general pull-back transformation (\ref{algtransf}) of a hypergeometric
equation is a Fuchsian equation. We are looking for situations when the
transformed equation is hypergeometric as well. In this Section we rather
look for transformed equations with at most 3 singular points. Since any
such Fuchsian equation can be transformed to a hypergeometric equation by
fractional-linear transformations, it is appropriate to ignore exact
location of singular points for a while. We loosely follow the 5-step
classification scheme in \cite[Section 3]{algtgauss}, with $N=3$, etc.

The requirement that the transformed equation must have at most 3 singular
points is restrictive. The covering $\varphi:\PP^1\to\PP^1$ essentially
determines singularities and local exponent differences of the transformed
equation. Here are basic general facts which we use (or refer to).
\begin{lemma} \label{genrami}
Let $\varphi:\PP^1_x\to\PP^1_z$ denote a finite covering of a projective
line $\PP^1_z$ with the rational parameter $z$ by a projective line
$\PP^1_x$ with the rational parameter $x$. Let $H_1$ denote a hypergeometric
equation on $\PP^1_z$, and let $H_2$ denote the pull-back transformation of
$H_1$ under $(\ref{algtransf})$. Let $d$ denote the degree of $\varphi$, and
let $S\in\PP^1_x$, $Q\in\PP^1_z$ be points such that $\varphi(S)=Q$.
\begin{enumerate}
\item If the point $Q$ is non-singular for $H_1$, then the point $S$ is
non-singular for $H_2$ only if the covering $\varphi$ does not branch at
$S$.%
\item If the point $Q$ is a singular point for $H_1$, then the point $S$ is
non-singular for $H_2$ only if the local exponent difference at $Q$ is equal
to $1/n$, where $n$ is the branching order of $\varphi$ at $S$.%
\item Let $\singset$ denote a set of $3$ points on $\PP^1_z$. If all
branching points of $\varphi$ lie above $\singset$, then there are exactly
$d+2$ distinct points on $\PP^1_x$ above $\singset$. Otherwise there are
more
than $d+2$ distinct points above $\singset$.%
\item Suppose that the equations $H_1$ and $H_2$ are hypergeometric. Let
$e_1,e_2,e_3$ denote the local exponent differences for $H_1$, and let
$e'_1,e'_2,e'_3$ denote the local exponent differences for $H_2$. Suppose
that the local exponent differences are real positive numbers, and that
$e_1+e_2+e_3\neq 1$. Then
\begin{equation} \label{trareas}
d=\frac{1-e'_1-e'_2-e'_3}{1-e_1-e_2-e_3}.
\end{equation}
\end{enumerate}
\end{lemma}
\proof The first two statements are weaker formulations of parts 2, 3 of
\cite[Lemma 2.4]{algtgauss}. The third statement is part 1 of \cite[Lemma
2.5]{algtgauss}, and the last statement is a weaker formulation of part 2 of
\cite[Lemma 2.5]{algtgauss}; they are consequences of Hurwitz' formula
\cite[Corollary IV.2.4]{harts}.\qed

Here are restrictions on coverings $\varphi:\PP^1\to\PP^1$ and local
exponent differences for algebraic transformations of hyperbolic
hypergeometric functions. They are stronger versions of the constraints used
in Step 2 of the classification scheme in \cite[Section 3]{algtgauss}.
\begin{lemma} \label{boundlemma}
Let $k,\ell,m$ denote positive integers such that
\begin{equation} \label{basicineq}
\frac{1}{k}+\frac{1}{\ell}+\frac{1}{m}<1 \qquad \mbox{and} \qquad
k\le\ell\le m.
\end{equation}
Let $H$ denote hypergeometric equation $(\ref{hpgde})$ such that the local
exponent differences are equal to $1/k,1/\ell,1/m$. Suppose that pull-back
transformation $(\ref{algtransf})$ transforms $H$ to a hypergeometric
equation. Let $d$ denote the degree of the covering $\varphi:\PP^1\to\PP^1$.
\begin{enumerate}
\item The points $x=0,1,\infty$ are actual singularities of the transformed
equation, and they lie above the subset $\{0,1,\infty\}$ of the
$z$-projective line. The covering $\varphi$ branches only above this subset,
so $\varphi(x)$ is a Belyi function.%
\item The following equality holds:
\begin{equation} \label{basiceq}
d-\left\lfloor \frac{d}k\right\rfloor-\left\lfloor \frac{d}\ell
\right\rfloor-\left\lfloor\frac{d}m\right\rfloor=1.
\end{equation}
\item The following inequality holds:
\begin{equation} \label{boundd}
d\left(1-\frac{1}{k}-\frac{1}{\ell}-\frac{1}{m}\right)\le
1-\frac{3}{m}.
\end{equation}
\item If $m>d$ then $1/d+1/k+1/\ell\ge 1$.
\item If $m\le d$ then
\begin{equation} \label{boundm}
\left(1-\frac1k-\frac1\ell\right)m^2-2m+3\le 0 \qquad\mbox{and} \qquad
\frac23\le \frac1k+\frac1\ell<1.
\end{equation}
\end{enumerate}
\end{lemma}
\proof Let $\singset$ denote the subset $\{0,1,\infty\}$ of the
$z$-projective line. By part 2 of Lemma \ref{genrami}, there are at most
$\lfloor d/k\rfloor$, $\lfloor d/\ell\rfloor$, $\lfloor d/m\rfloor$
non-singular points above the $z$-points with the local exponent differences
$1/k$, $1/\ell$, $1/m$ respectively. By part 3 of Lemma \ref{genrami}, the
number of singular points above $\singset$ is at least
\begin{equation} \label{numspoints}
d+2-\left\lfloor \frac{d}k\right\rfloor-\left\lfloor \frac{d}\ell
\right\rfloor-\left\lfloor\frac{d}m\right\rfloor.
\end{equation}
This number is greater than $2+d\,(1-1/k-1/\ell-1/m)$, so it is at least
$3$. On the other hand, the transformed equation has at most $3$ singular
points. Therefore the transformed equation has exactly $3$ singular points,
expression (\ref{numspoints}) is equal to $3$, and $\varphi$ does not branch
outside $\singset$ by part 1 of Lemma \ref{genrami}. Parts 1 and 2 of this
Lemma follow.

We rewrite formula (\ref{basiceq}) as follows:
\begin{equation}
d\left( 1-\frac1k-\frac1\ell-\frac1m \right)+T=1,
\end{equation}
where $T$ is the sum of positive local exponent differences at the singular
points of the transformed equation. (This is equivalent to (\ref{trareas}),
with $T=e'_1+e'_2+e'_3$.) We have $T\ge 3/m$, which implies inequality
(\ref{boundd}).

If $m>d$ then we use formula (\ref{basiceq}) to derive
\begin{equation}
1=d-\left\lfloor \frac{d}k\right\rfloor-\left\lfloor
\frac{d}\ell\right\rfloor\ge d\left(1-\frac1k-\frac1\ell\right),
\end{equation}
which gives part 4 of this Lemma. If $m\ge d$, we derive the first
inequality in (\ref{boundm}) after replacing $d$ by $m$ in (\ref{boundd}).
We have $1/k+1/\ell<1$ from (\ref{basicineq}). The quadratic expression in
$m$ achieves non-positive values only if $1/k+1/\ell\ge 2/3$ (consider the
discriminant). \qed

These restrictions essentially give a finite list of possibilities for the
integer tuple $(k,\ell,m,d)$. Indeed, inequality (\ref{boundd}) bounds $d$
once $k,\ell,m$ are fixed, and part 5 of Lemma \ref{boundlemma} gives
finitely many possibilities for the triple $(k,\ell,m)$. Only when $m>d$ we
formally have infinitely many possibilities; but then we expect to arrive at
specializations of algebraic transformations with unrestricted parameters.
The inequalities give the following possibilities:
\begin{eqnarray*}
(2,\ell,m,2),\quad (2,3,m,3...6),\quad (2,4,m,4),\quad (3,3,m,3),\quad(2,3,7,7...24),\\
\hspace{-11pt}(2,3,8,8...15),\quad(2,3,9,9...12),\quad(2,3,10,10),\quad(2,4,5,5...8),\quad(2,4,6,6).
\end{eqnarray*}
Here $d$ is sometimes represented by an integer interval of possible values,
and the unevaluated parameters $\ell,m$ can be large enough integers.
Formula (\ref{basiceq}) rejects some of these possibilities.

The next step is to produce a list of possible branching patterns. Because
of parts 1 and 2 of Lemma \ref{boundlemma}, we have to take the maximal
possible number $\lfloor d/k\rfloor$, $\lfloor d/\ell\rfloor$ or $\lfloor
d/m\rfloor$ of non-singular points above the 3 singular $z$-points. The
remaining residual branches above $z=0,1,\infty$ should coalesce into
precisely 3 distinct points. In particular, we have to ignore the cases when
there remains less than 3 residual branches. The final list of possible
branching patterns is presented in the first three columns of Table
\ref{figtabc} (for the cases with $m>d$) and Table \ref{figtab} (for the
cases with $m\le d$). In Table \ref{figtabc} we ignore fractional-linear
transformations, and we drop the condition $\ell\le m$ for degree 2
transformations.
\begin{table}
\begin{center} \begin{tabular}{|c|c|c|c|c|}
\hline \multicolumn{2}{|c|}{Local exponent differences} & Degree & Covering
& Coxeter \\ \cline{1-2}
$(1/k,\,1/\ell,\,1/m)$ & above & $d$ & composition & decomposition \\
 \hline\hline
$(1/2,\,1/\ell,\,1/m)$ & $(1/\ell,\,1/\ell,\,2/m)$ & 2 & indecomposable & yes \\
$(1/2,\,1/3,\,1/m)$ & $(1/2,\,1/m,\,2/m)$ & 3 & indecomposable & yes \\
$(1/2,\,1/3,\,1/m)$ & $(1/3,\,1/m,\,3/m)$ & 4 & indecomposable & yes \\
$(1/2,\,1/3,\,1/m)$ & $(1/3,\,2/m,\,2/m)$ & 4 & \multicolumn{2}{c|}{no covering}\\
$(1/2,\,1/3,\,1/m)$ & $(1/m,\,1/m,\,4/m)$ & 6 & $2\times 3$ & yes \\
$(1/2,\,1/3,\,1/m)$ & $(2/m,\,2/m,\,2/m)$ & 6 & $2\times 3$ or $3\times 2$ & yes \\
$(1/2,\,1/3,\,1/m)$ & $(1/m,\,2/m,\,3/m)$ & 6 & \multicolumn{2}{c|}{no covering}\\
$(1/2,\,1/4,\,1/m)$ & $(1/m,\,1/m,\,2/m)$ & 4 & $2\times 2$ & yes \\
$(1/3,\,1/3,\,1/m)$ & $(1/m,\,1/m,\,1/m)$ & 3 & indecomposable & no \\
\hline
\end{tabular} \end{center}
\caption{Classical transformations of hyperbolic hypergeometric functions}
\label{figtabc}
\end{table}
Branching patterns are uniquely determined by the starting local exponent
differences $(1/k,1/\ell,1/m)$, transformed local exponent differences, and
the stated principle to have maximal number of non-singular points above
$z=0,1,\infty$. For example, branching pattern for the degree 4
transformation of Table \ref{figtabc} between hypergeometric equations with
the local exponent differences $(1/2,\,1/3,\,1/m)$ and $(1/3,\,1/m,\,3/m)$
can be schematically denoted by $2+2=3+1=3+1$. This means that all points
above the $z$-point with the local exponent difference $1/2$ have branching
order 2, and that there must be a branching point with order 3 and a
non-branching point above each of the other two points. For one more
example, the degree 9 covering of Table \ref{figtab} has the branching
pattern $2+2+2+2+1=3+3+3=7+1+1$. The same notation for branching pattern is
used in \cite{algtgauss}.
\begin{table}
\begin{center} \begin{tabular}{|c|c|c|c|c|}
\hline \multicolumn{2}{|c|}{Local exponent differences} & Degree & Covering
& Coxeter \\ \cline{1-2}
$(1/k,\,1/\ell,\,1/m)$ & above & $d$ & composition & decomposition \\
 \hline\hline
$(1/2,\,1/3,\,1/7)$ & $(1/3,\,1/3,\,1/7)$ & 8 & indecomposable & no \\
$(1/2,\,1/3,\,1/7)$ & $(1/2,\,1/7,\,1/7)$ & 9 & indecomposable & no \\
$(1/2,\,1/3,\,1/7)$ & $(1/3,\,1/7,\,2/7)$ & 10 & indecomposable & yes \\
$(1/2,\,1/3,\,1/7)$ & $(1/7,\,1/7,\,3/7)$&12&\multicolumn{2}{c|}{no covering}\\
$(1/2,\,1/3,\,1/7)$ & $(1/7,\,2/7,\,2/7)$&12&\multicolumn{2}{c|}{no covering}\\
$(1/2,\,1/3,\,1/7)$ & $(1/3,\,1/7,\,1/7)$&16&\multicolumn{2}{c|}{no covering}\\
$(1/2,\,1/3,\,1/7)$ & $(1/7,\,1/7,\,2/7)$ & 18 & $2\times 9$ & no \\
$(1/2,\,1/3,\,1/7)$ & $(1/7,\,1/7,\,1/7)$ & 24 & $3\times 8$ & yes \\
$(1/2,\,1/3,\,1/8)$ & $(1/3,\,1/8,\,1/8)$ & 10 & indecomposable & no \\
$(1/2,\,1/3,\,1/8)$ & $(1/4,\,1/8,\,1/8)$ & 12 & $2\times 2\times 3$ & yes \\
$(1/2,\,1/3,\,1/9)$ & $(1/9,\,1/9,\,1/9)$ & 12 & $3\times 4$ & no \\
$(1/2,\,1/4,\,1/5)$ & $(1/4,\,1/4,\,1/5)$ & 6 & indecomposable & no \\
$(1/2,\,1/4,\,1/5)$ & $(1/5,\,1/5,\,1/5)$ & 8 &\multicolumn{2}{c|}{no
covering}\\ \hline
\end{tabular} \end{center}
\caption{Non-classical transformations of hyperbolic hypergeometric
functions} \label{figtab}
\end{table}

Now we have to determine all coverings which have those branching patterns.
Given a branching pattern, there is often exactly one covering with that
branching pattern up to fractional-linear transformations. But not for any
branching pattern a covering exists, and there can be several different
coverings with the same branching pattern. Section \ref{fincovers} is
devoted to computing coverings with a given branching pattern. First we
outline there a straightforward method with undetermined coefficients, which
is feasible if $d\le 6$. Then we introduce a more appropriate algorithm,
which was actually used (within computer algebra system {\sf Maple}) to
compute coverings for Table \ref{figtab}. Final information about existing
coverings is given in the fourth columns of Table \ref{figtabc} and Table
\ref{figtab}. It turns out that for any candidate branching pattern there is
at most one covering up to fractional-linear transformations. If $m>d$ (see
Table \ref{figtabc}), we get the coverings of the classical algebraic
transformations due to Gauss, Euler, Kummer, Goursat. If $m\le d$ (see Table
\ref{figtab}), we get new coverings of degree 6, 8, 9, 10, 12, 18, 24.
Existence of some of these coverings is shown in
\cite{hodgkins2,beukers,kitaevdb}. For both Tables, it was straightforward
to figure out possible compositions of small degree coverings and identify
them with the unique coverings for suitable branching patterns. Numbers in
the multiplicative notation for decomposable coverings mean degrees of
constituent coverings, as in \cite{algtgauss}.

The last step is to determine algebraic transformations of hypergeometric
functions with the rational argument determined by a computed covering. The
factor $\theta(x)$ in (\ref{algtransf}) should shift local exponents at
potentially non-singular points to the values 0 and 1, and it should shift
one local exponent at both $x=0$ and $x=1$ to the value 0. A suitable
pull-back transformation induces a hypergeometric identity like
(\ref{hpgtransf}) for each singular $x$-point $S$ which lies above a
singular $z$-point. To achieve this, one has to move the points $S$ and
$\varphi(S)$ to the locations $x=0$ and $z=0$ respectively (by
fractional-linear transformations), and identify the two solutions with the
local exponent 0 and the value 1 at \mbox{$x=0$} and $z=0$ respectively. It
is convenient to use Riemann's $P$-notation for these purposes; see
\cite[Section 3.9]{specfaar} or \cite[Section 2]{algtgauss}. Each
positioning of $x=0$ above $z=0$ gives a few hypergeometric identities like
(\ref{hpgtransf}). First of all, we have Euler's and Pfaff's
fractional-linear transformations \cite[Theorem 2.2.5]{specfaar}, which
permute other two singular points and their local exponents. Additionally,
simultaneous permutation of the local exponents at $x=0$ and $z=0$ gives the
following hypergeometric identity.
\begin{lemma} \label{basislem}
Suppose that a pull-back transformation induces identity $(\ref{hpgtransf})$
in an open neighborhood of $x=0$. Then $\varphi(x)^{1-C}\sim K
x^{1-\widetilde{C}}$ as $x\to 0$ for some constant $K$, and the following
identity holds (if both hypergeometric functions are well-defined):
\[ 
\hpg{2}{1}{\!1+\widetilde{A}-\widetilde{C},1+\widetilde{B}-\widetilde{C}}
{2-\widetilde{C}}{\,x}=\theta(x)\frac{\varphi(x)^{1-C}}{K\,x^{1-\widetilde{C}}}
\,\hpg{2}{1}{\!1+A-C,1+B-C}{2-C}{\varphi(x)}.
\] 
\end{lemma}
\proof This is Lemma 2.3 in \cite{algtgauss}. \qed

As it turns out, algebraic transformations for Table \ref{figtabc} (i.e.,
the case $m>d$) are special cases of the classical transformations due to
Gauss, Euler, Kummer, Goursat. We give a few instances of these
transformations in Section \ref{hyperids}. Algebraic transformations for
Table \ref{figtab} (i.e., the case $m\le d$) are modern, though some of them
are predicted in \cite{hodgkins2}. We present these transformations (up to
Euler's and Pfaff's fractional-linear transformations, and Lemma
\ref{basislem}) in Section \ref{hyperids} as well.

The rest of this Section is devoted to explaining the last columns of Tables
\ref{figtabc} and \ref{figtab}. Recall \cite{yoshida,beukers} that a {\em
Schwarz map} for a hypergeometric equation $H_1$ is an analytic map from
the upper half-plane \mbox{$\HH=\{z\in\CC\,|\,\mbox{Im }z>0\}$} given by a
quotient of two solutions of $H_1$. If the local exponent differences
$e_0,e_1,e_\infty$ of $H_1$ are real numbers in the interval $[0,1)$, then
the image of a Schwarz map is a curvilinear triangle on the Riemann sphere.
Such a triangle is called {\em Schwarz triangle}.
The vertices are images of the 3 singular points, and the angles there are
equal to $\pi e_0$, $\pi e_1$, $\pi e_\infty$ correspondingly; the sides are
circular arcs. Analytic continuation of a Schwarz map follows the {\em
Schwarz reflection principle}: the image of the other half-plane under
analytic continuation across $(0,1)$, $(1,\infty)$ or $(-\infty,0)$ is a
fractional-linear reflection of the Schwarz triangle across the
corresponding side of itself.

In our case, the local exponent differences are $1/k,1/\ell,1/m$, and
$1/k+1/\ell+1/m<1$. The sides of a Schwarz triangle are geodesic curves
with respect to a hyperbolic metric on the Riemann sphere, defined on some
Poincare disk. 
Repeated analytic continuation gives a tessellation of the Poincare disk
into curvilinear triangles with the angles $\pi/k,\pi/\ell,\pi/m$.

\begin{figure}
\[
\setlength{\unitlength}{0.75pt}
\begin{picture}(600,306)
\put(-60,18){\resizebox{10.2cm}{!}{\includegraphics{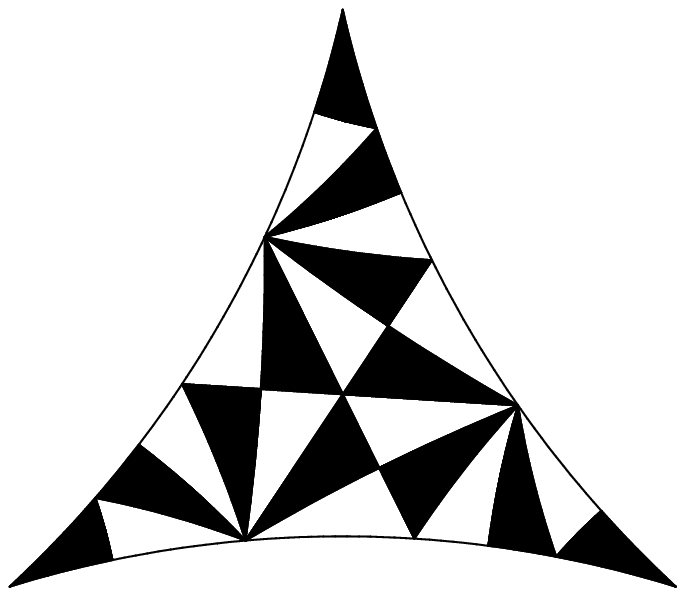}}}
\put(28,128){\resizebox{5.4cm}{!}{\includegraphics{schwartz237a.eps}}}
\put(4,-60){\resizebox{6cm}{!}{\includegraphics{schwartz238a.eps}}}
\put(126,30){\it(a)} \put(478,160){\it(b)} \put(478,-6){\it(c)}
\end{picture}
\]
\caption{Coxeter decompositions of hyperbolic triangles}
\label{triangles237}
\end{figure}

Consider a pull-back transformation of the hypergeometric equation $H_1$ to
a hypergeometric equation $H_2$, of degree $d$. Suppose that its covering
\mbox{$\varphi:\PP^1\to\PP^1$} is defined over $\RR$, and that it branches
only above the singular points of $H_1$. If $s:\HH\to\CC$ is a Schwarz map
for $H_2$, then a branch of $s\circ\varphi^{-1}$ is a Schwarz map for $H_1$.
The Schwarz triangles of the $d$ branches of $s\circ\varphi^{-1}$ tessellate
the Schwarz triangle of $s$, like in Figure \ref{triangles237}. In this case
the degree expression (\ref{trareas}) can be interpreted as the quotient of
areas of Schwarz triangles for the two hypergeometric equations, in the
hyperbolic or spherical metric. Transformations of hypergeometric equations
which admit these tessellations are implicitly classified in
\cite{hodgkins2,beukers}. Tessellations of hyperbolic triangles and
quadrangles into hyperbolic triangles are classified in
\cite{felixon,nyhyperb}, where they are called {\em Coxeter decompositions}
and {\it divisible tilings} respectively. The classification in
\cite{nyhyperb} is incomplete; for instance, it misses triangulation
{\em(b)} in Figure \ref{triangles237}. We adopt the terminology of
\cite{felixon}.

The last columns of Tables \ref{figtabc} and \ref{figtab} tell us which
transformations of hypergeometric equations with hyperbolic solutions admit
Coxeter decompositions of Schwarz triangles. In particular, the three such
transformations in Table \ref{figtab} are anticipated in
\cite{hodgkins2,beukers}. Their Coxeter decompositions are depicted in
Figure \ref{triangles237}. All classical transformations except one cubic
transformation admit these tessellations.

\section{Computation of Belyi functions}
\label{fincovers}

Here we consider the problem of computing finite coverings
$\varphi:\PP^1\to\PP^1$ of given degree $d$ and with a given branching
pattern. We assume that all branching points lie above a set $\singset$ of 3
points, so the desired functions $\varphi(x)$ are Belyi functions. By part 3
of Lemma \ref{genrami}, there must be exactly $d+2$ distinct points above
$\singset$. We expect finitely many (or no) solutions to this problem.
 Algorithm \ref{compcovering} of this Section was used to
compute transformations implied by Tables \ref{figtabc} and \ref{figtab}. As
mentioned, there is at most one solution for branching patterns there.

First we outline a naive method with undetermined coefficients, which is
feasible if $d\le 6$. To fix ideas, consider the branching pattern
$2+2+2=3+3=2+2+2$ for the sixth entry 
in Table \ref{figtabc}. Up to fractional-linear transformations of the
$z$-line, we may assume that points with these branching orders lie above
$z=1,0$ and $\infty$ respectively. We choose the 
points above $z=\infty$ to be $x=0$, $x=1$, $x=\infty$. Then the Belyi
function should have the form
\begin{equation}  \label{phi6}
\varphi(x)=
\frac{(u_2x^2+u_1x+u_0)^3}{x^2(x-1)^2},
\end{equation}
where $u_2,u_1,u_0$ are undetermined, and the roots of $u_2x^2+u_1x+u_0$ are
the points above $z=0$. 
The branching pattern above $z=1$ implies that the numerator of
$\varphi(x)-1$ must be a square of a cubic 
polynomial $P(x)$. This condition gives 7 polynomial equations in the 4
coefficients of $P(x)$ and $u_2,u_1,u_0$. These equations can be feasibly
solved with assistance of a computer algebra package if $d$ is not large.

To compute Belyi functions more efficiently, we propose to pull-back the
differential $dz/z$ with respect to $\varphi$, still with undetermined
coefficients as in (\ref{phi6}). The poles of the pull-backed differential
are simple, and they are located at the points above $z=0$ and $z=\infty$.
Its zeroes are the branching points which do not lie above $z=0$ or
$z=\infty$; the multiplicities of the zeroes are 1 less than the
corresponding branching orders. In our example, all those branching points
must lie above $z=1$, so they are the roots of $P(x)$. Moreover, they must
be simple roots of $P(x)$, so we get this polynomial just by computing the
pull-back of $dz/z$. Explicitly, the pull-back of $dz/z$ is equal to
\begin{equation} \label{logdif1}
\frac{\varphi'(x)}{\varphi(x)}\,dx=
\frac{2u_2x^3-(4u_2\!+\!u_1)\,x^2-(u_1\!+\!4u_0)\,x+2u_0}
{x\,(x-1)\,(u_2x^2+u_1x+u_0)}\,dx.
\end{equation}
Let $\widetilde{P}(x)$ denote the polynomial in the numerator on the
right-hand side; it must be proportional to $P(x)$. Therefore $\varphi(x)-1$
is proportional to $\widetilde{\varphi}(x):=\widetilde{P}(x)^2/x^2(x-1)^2$.
Further, consider the pull-back of the differential \mbox{$d(z-1)/(z-1)$}:
\begin{equation} \label{logdif2}
\frac{\widetilde{\varphi}\,'(x)}{\widetilde{\varphi}(x)}\,dx=4\,\frac{
u_2x^4-2u_2x^3+(2u_2\!+\!u_1\!+\!2u_0)x^2-2u_0x+u_0}
{x\,(x-1)\,(2u_2x^3-(4u_2\!+\!u_1)x^2-(u_1\!+\!4u_0)\,x+2u_0)}\,dx.
\end{equation}
By the same reasoning, the zeroes of this differential are the branching
points of $\varphi$ which do not lie above $z=1$ or $z=\infty$, with the
multiplicities diminished by 1. In our case, those branching points must lie
above $z=0$. Hence the polynomial in the numerator of (\ref{logdif2}) is
proportional to $(u_2x^2+u_1x+u_0)^2$. This gives easy polynomial equations
in $u_2$, $u_1$, $u_0$. Since we want $u_2u_0\neq 0$,
\begin{eqnarray} \label{exshows}
\frac{u_2^2}{u_2}=-\frac{2u_1u_2}{2u_2}=\frac{u_1^2+2u_0u_2}
{2u_2\!+\!u_1\!+\!2u_0}=-\frac{2u_0u_1}{2u_0}=\frac{u_0^2}{u_0}.
\end{eqnarray}
We solve that $u_2=-u_1=u_0$. Therefore $\varphi(x)$ is proportional to
\mbox{$(x^2\!-x+1)^3/x^2(x-1)^2$}. The scalar multiple can be found from the
condition that $z=1$ is the third branching locus of $\varphi(x)$. We
derive:
\[ 
\varphi(x)=\frac{4}{27}\frac{(x^2-x+1)^3}{x^2\,(x-1)^2}, \qquad
\varphi(x)-1=\frac{(x+1)^2\,(2x-1)^2\,(x-2)^2}{27\,x^2\,(x-1)^2}.
\] 
We have solved the problem by hand! The solution is unique up to
fractional-linear transformations. Note that there are two different ways to
compose coverings of degree 2 and 3 and get a covering with the considered
branching pattern; see Table \ref{figtabc} and \cite[Section 4]{algtgauss}.
Up to fractional-linear transformations, those two compositions must give
the same covering computed here.

Now we present general Algorithm \ref{compcovering} for finding Belyi
functions with a given branching pattern. To formulate it more conveniently,
we restrict ourselves to branching patterns that are relevant for the
purposes of this paper. Note that for all transformations for Table
\ref{figtab} (and almost all transformations for Table \ref{figtabc}) there
is a $z$-point with the local exponent difference $1/2$. For coverings of
these transformations assumption {\em(b)} of Algorithm \ref{compcovering}
holds. If this assumption is dropped, then Step 1 should try to assign the
fiber with smallest branching orders to $z=1$, the function
$\widetilde{\varphi}(x)$ in Step 3 has a more complicated form, and more
undetermined coefficients are needed.
\begin{alg} \label{compcovering} \rm
{\em Input:} a branching pattern (that is, 3 collections of branching
orders) and degree $d$. We assume:\\
{\em (a)} the branching orders in the same fiber sum up to $d$, and
there are $d+2$ branching orders in total;\\
{\em (b)} one of the 3 collections prescribes only branching orders $2$
and at most one unramified (i.e., simple, not branching) point.\\
{\em Output:} All Belyi functions (up to fractional-linear transformations)
whose coverings branch only above the set $\singset=\{0,1,\infty\}$
with the given branching orders.\\
{\bf Step 1.} Prescribe the branching orders mentioned in assumption {\em
(b)} to the fiber of the point $z=1$, and prescribe other two collections to
the fibers of $z=0$ and $z=\infty$. Choose the points $x=0$, $x=1$,
$x=\infty$ above $\singset$ in a convenient way: if an unramified point is
prescribed above $z=1$, choose it to be $x=\infty$; see also remarks
immediately below. Consider the other points above $z=0$ and $z=\infty$ as
unknown. Accordingly, write $\varphi(x)=K\,P(x)/Q(x)$, where $K$ is an
undetermined constant, and $P(x)$, $Q(x)$ are monic polynomials in the
square-free factorized form (following the branching pattern) with some
undetermined coefficients.\\
{\bf Step 2.} Compute the pull-back $\varphi'(x)\,dx/\varphi(x)$ of $dz/z$.
Let $R(x)$ be the numerator of $\varphi'(x)/\varphi(x)$. The roots of $R(x)$
are the branching points above $z=1$.\\
{\bf Step 3.} Let $\widetilde{\varphi}(x)=R(x)^2/Q(x)$ and compute
the rational function
$\Phi(x)=\widetilde{\varphi}{\,'}(x)/\widetilde{\varphi}(x)$. The
numerator of $\Phi(x)$ has the same roots as the polynomial
$P(x)$, but their multiplicity is 1 less than in $P(x)$. This
gives a set of algebraic
equations in the undetermined coefficients.\hspace{1pt}\\
{\bf Step 4.} Solve the algebraic equations obtained in Step 3 by Gr\"obner
basis methods, and find possible pairs of polynomials $P(x),Q(x)$ with the
right factorization pattern. For each non-degenerate solution, the constant
$K$ in the target $\varphi(x)=K\,P(x)/Q(x)$ is such that the function
$1-K\,P(x)/Q(x)$ is proportional to $\widetilde{\varphi}(x)$. If necessary,
compose the output functions $\varphi(x)$ with suitable fractional-linear
transformations to move some $x$-points (or even $z$-points) to final
desired locations.\hfill$\Box$
\end{alg}
When applying this algorithm to the entries of Table \ref{figtab}, it is
convenient to choose the points $x=0$, $x=1$, $x=\infty$ in Step 1 to be the
singular points of the transformed hypergeometric equation. In general, a
good strategy for Step 1 is to choose points which have different branching
orders than the most points in the same fiber. On the other hand, the
algorithm can be more effective if we choose points with maximal branching
orders as $x=0$, $x=1$, $x=\infty$. With this modification, the function
$\widetilde{\varphi}(x)$ in Step 3 may acquire an extra linear factor. It
may be convenient not to make a choice for $x=1$. Then we would have an
extra variable and the algebraic equations would be weighted-homogeneous
(respecting to the transformations $x\mapsto\alpha x$). The extra degree of
freedom can be used to avoid complicated algebraic numbers.

We note that Step 3 produces enough algebraic equations between the
undetermined coefficients, because the restrictions on the polynomials
$P(x),Q(x),R(x)$ and the denominator of $\Phi(x)$ determine the desired
branching pattern for $\varphi(x)$. Hence the algorithm is correct. As
example (\ref{exshows}) shows, the set of equations is  likely to be
overdetermined, which only helps in Gr\"obner basis computations.

Compared with the naive method described at the beginning of this Section,
algebraic equations of Algorithm \ref{compcovering} have fewer undetermined
coefficients, lower degree, and fewer degenerate (or {\em parasytic}
\cite{kreines}) solutions. The computations are still tedious, but all
coverings of Table \ref{figtab} were computed using the computer algebra
package {\sf Maple} in a matter of hours. There is an article
\cite{shabatgb} where quadratic differentials are used to characterize some
Belyi maps. But \cite{kreines,kitaevdb} exploit the naive method.

\section{Hypergeometric identities}
\label{hyperids}

Here we present our main results. We give all algebraic transformations for
Table \ref{figtab}, up to Euler's and Pfaff's fractional-linear
transformations and Lemma \ref{basislem}. But first we exhibit a few
classical transformations.

Relevant instances of quadratic transformations can be obtained by setting
$a=1/2-1/\ell-1/m$, $b=1/2-1/\ell+1/m$ in formulas
(\ref{hpgtr2})--(\ref{hpgtr2a}). Examples of classical transformations of
degree 3 or 4 are:
\begin{eqnarray} \label{fcomp6a}
\hpg{2}{1}{a,\,\frac{1-a}{3}\,}{\frac{4a+5}{6}}{\,x} & \equal &
\left(1-4x\right)^{-a}\,
\hpg{2}{1}{\frac{a}{3},\,\frac{a+1}{3}}{\frac{4a+5}{6}}
{\frac{27\,x}{(4x-1)^3}},\\
\label{fcomp8a}\hpg{2}{1}{\frac{4a}{3},\,\frac{4a+1}{3}}{\frac{4a+5}{6}}
{\,x} & \equal & \left(1+8x\right)^{-a}\,
\hpg{2}{1}{\frac{a}{3},\,\frac{a+1}{3}}{\frac{4a+5}{6}}
{\frac{64\,x\,(1-x)^3}{(1+8x)^3}},\\
\label{cubic3} \hpg{2}{1}{c,\;\frac{c+1}{3}}{\frac{2c+2}{3}}{\,x}
& \equal & \left(1+\omega^2 x\right)^{-c}\,
\hpg{2}{1}{\frac{c}{3},\,\frac{c+1}{3}}{\frac{2c+2}{3}}
{\frac{3(2\omega\!+\!1)\;x(x-1)}{(x+\omega)^3}}.
\end{eqnarray}
Here $\omega$ is a primitive cubic root of unity (so $\omega^2+\omega+1=0$),
$a=1/4\pm 3/2m$ and $c=1/2\pm 3/2m$. These identities correspond to the
indecomposable pull-back coverings $\varphi:\PP^1\to\PP^1$ of Table
\ref{figtabc}. For more complete lists of classical algebraic
transformations we refer to \cite{goursat,algtgauss}.

Now we present non-classical transformations of hyperbolic hypergeometric
functions. Our list of transformations is basically complete, it was
computed following the plan of Section \ref{otherat}. If a covering for
Table \ref{figtab} is indecomposable, all corresponding two-term
hypergeometric identities can be obtained from the exhibited below, by using
Euler's and Pfaff's fractional-linear transformations and Lemma
\ref{basislem}. If the covering is decomposable, we indicate a composition
of hypergeometric identities of smaller degree.

A covering for degree 8 pull-back transformations between hypergeometric
equations with the local exponent differences $(1/2,\,1/3,\,1/7)$ and
$(1/3,\,1/3,\,1/7)$ is given by
\begin{equation}
\varphi_1(x)=\frac{x\,(x-1)\,\left(27x^2-(723\!+\!1392\omega)x-496\!+\!696\omega\right)^3}
{64\,\big((6\omega+3)x-8-3\omega\big)^7}.
\end{equation}
Here $\omega$ satisfies $\omega^2+\omega+1=0$ as in formula (\ref{cubic3}).
Note that the conjugation $\omega=-1-\omega$ acts in the same way as a
composition with fractional-linear transformation interchanging the points
$x=0$ and $x=1$. This confirms uniqueness of the covering. The covering is
computed in \cite{kitaevdb} as well. Here are hypergeometric identities:
\begin{eqnarray}
\hpg{2}{1}{2/21,\,5/21}{2/3}{\,x} & \equal &
{\textstyle\left(1-\frac{33+39\omega}{49}\,x\right)}^{-1/12}
\,\hpg{2}{1}{1/84,\,13/84}{2/3}{\varphi_1(x)},\\ \label{transf78}
\hpg{2}{1}{2/21,\,3/7}{6/7}{\,x} & \equal & (1-x)^{-1/84}
{\textstyle\left(1-\frac{241+464\omega}{9}\,x-\frac{8(62-87\omega)}
{27}\,x^2\right)}^{-1/28} \nonumber\\
& & \times\,\hpg{2}{1}{1/84,\,29/84}{6/7}{\frac1{\varphi_1(1/x)}}.
\end{eqnarray}

A covering for degree 9 pull-back transformations between hypergeometric
equations with the local exponent differences $(1/2,\,1/3,\,1/7)$ and
$(1/2,\,1/7,\,1/7)$ is given by
\begin{equation}
\varphi_2(x)=\frac{27\,x\,(x-1)\,(49x-31-13\xi)^7}
{49\,(7203x^3+(9947\xi-5831)x^2-(9947\xi+2009)x+275-87\xi)^3}.
\end{equation}
Here $\xi$ satisfies $\xi^2+\xi+2=0$. Hypergeometric identities
are:
\begin{eqnarray} \label{transf79}
\hpg{2}{1}{3/28,\,17/28}{6/7}{\,x} & \equal & \textstyle{\left(
1+\frac{7(10-29\xi)}{8}\,x-\frac{343(50-29\xi)}{512}\,x^2
+\frac{1029(362+87\xi)}{16384}\,x^3\right)}^{-1/28}\nonumber\\
&& \times\,\hpg{2}{1}{1/84,\,29/84}{6/7}{\varphi_2(x)},\\
\hpg{2}{1}{3/28,\,1/4}{1/2}{\,x} & \equal & \textstyle{\left(
1-\frac{17-29\xi}{21}\,x-\frac{41+203\xi}{147}\,x^2+
\frac{275-87\xi}{7203}\,x^3\right)}^{-1/28}\nonumber\\
&& \times\,\hpg{2}{1}{1/84,\,29/84}{1/2}{1-\varphi_2(1/x)}.
\end{eqnarray}

A covering for degree 10 pull-back transformations between hypergeometric
equations with the local exponent differences $(1/2,\,1/3,\,1/7)$ and
$(1/3,\,1/7,\,2/7)$ is given by
\begin{equation}
\varphi_3(x)=-\frac{x^2\,(x-1)\,(49x-81)^7}{4\,(16807x^3-9261x^2-13851x+6561)^3}.
\end{equation}
Hypergeometric identities are:
\begin{eqnarray}
\hpg{2}{1}{\!5/42,\,19/42}{5/7}{\,x} & \!\equal\! &
{\textstyle\left(1\!-\!\frac{19}9x\!-\!\frac{343}{243}x^2\!+\!\frac{16807}{6561}x^3\right)}^{-1/28}
\hpg{2}{1}{\!1/84,\,29/84}{6/7}{\varphi_3(x)}\!,\\
\hpg{2}{1}{\!5/42,\,19/42}{6/7}{\,x} & \!\equal\! &
{\textstyle\left(1\!-\!\frac{141}2x\!+\!\frac{5145}{32}x^2\!-\!\frac{16807}{256}x^3\right)}^{-1/28}
\hpg{2}{1}{\!1/84,29/84}{6/7}{\varphi_3(1\!-\!x)}\!,\nonumber\\ \\
\hpg{2}{1}{\!5/42,\,17/42}{2/3}{\,x} & \!\equal\! &
(1-x)^{-1/84}{\textstyle\left(1-\frac{81}{49}x\right)}^{-1/12}\,
\hpg{2}{1}{1/84,\,13/84}{2/3}{\frac1{\varphi_3(1/x)}}.
\end{eqnarray}

Degree 18 transformations between hypergeometric equations with the local
exponent differences $(1/2,\,1/3,\,1/7)$ and $(1/7,\,1/7,\,2/7)$ are
compositions of degree 9 and degree 2 transformations. The intermediate
hypergeometric equation has the local exponent differences $(1/2,1/7,1/7)$.
To get a hypergeometric identity, one can compose formula (\ref{hpgtr2a})
with $a=3/14$, $b=1/2$ and formula (\ref{transf79}).

Degree 24 transformations between hypergeometric equations with the local
exponent differences $(1/2,\,1/3,\,1/7)$ and $(1/7,\,1/7,\,1/7)$ are
compositions of degree 8 and degree 3 transformations. The intermediate
hypergeometric equation has the local exponent differences $(1/3,1/3,1/7)$.
Note that we have here a composition of two pull-back transformations which
do not admit a Coxeter decomposition, but the composite transformation does
admit a Coxeter decomposition. To get a hypergeometric identity, one can
compose formula (\ref{cubic3}) with $c=2/7$ and formula (\ref{transf78});
see \cite[formula (76)]{algtgauss}.

A covering for degree 10 pull-back transformations between hypergeometric
equations with the local exponent differences $(1/2,\,1/3,\,1/8)$ and
$(1/3,\,1/8,\,1/8)$ is given by
\begin{equation}
\varphi_4(x)=\frac{4\,x\,(x-1)\,(8\beta x+7-4\beta)^8} {(2048\beta
x^3-3072\beta x^2-3264x^2+912\beta x+3264x+56\beta-17)^3}.
\end{equation}
Here $\beta$ satisfies $\beta^2+2=0$. Hypergeometric identities
are:
\begin{eqnarray}
\hpg{2}{1}{\!5/24,13/24}{7/8}{\,x} &\equal &
{\textstyle\left(1+\frac{16(4-17\beta)}9x-\frac{64(167-136\beta)}{243}x^2
+\frac{2048(112-17\beta)}{6561}x^3\right)}^{-1/16}\nonumber\\
&&\times\,\hpg{2}{1}{1/48,\,17/48}{7/8}{\varphi_4(x)},\\
\hpg{2}{1}{5/24,\,1/3}{2/3}{\,x} &\equal & (1-x)^{-1/48}
{\textstyle\left(1-\frac{8+7\beta}{16}x\right)}^{-1/6}
\hpg{2}{1}{\!1/48,7/48}{2/3}{\frac1{\varphi_4(1/x)}}.\nonumber\\
\end{eqnarray}

Degree 12 transformations between hypergeometric equations with the local
exponent differences $(1/2,\,1/3,\,1/8)$ and $(1/4,\,1/8,\,1/8)$ are
compositions of a degree 3 transformation and two quadratic transformations.
The intermediate hypergeometric equations have the local exponent
differences $(1/2,1/4,1/8)$ and $(1/2,1/8,1/8)$. To get a hypergeometric
identity, one can compose formula (\ref{hpgtr2a}) with $a=1/4$, $b=1/2$,
formula (\ref{hpgtr2}) with $a=1/8$, $b=5/8$, and formula (\ref{fcomp6a})
with $a=1/16$.

Degree 12 transformations between hypergeometric equations with the local
exponent differences $(1/2,\,1/3,\,1/9)$ and $(1/9,\,1/9,\,1/9)$ are
compositions of degree 4 and degree 3 transformations. The intermediate
hypergeometric equation has the local exponent differences $(1/3,1/3,1/9)$.
To get a hypergeometric identity, one can compose formula (\ref{cubic3})
with $c=1/3$ and formula (\ref{fcomp8a}) with $a=1/12$.

A covering for degree 6 pull-back transformations between hypergeometric
equations with the local exponent differences $(1/2,\,1/4,\,1/5)$ and
$(1/4,\,1/4,\,1/5)$ is given by
\begin{equation}
\varphi_5(x)=\frac{4i\,x\,(x-1)\,(4x-2-11i)^4}{(8x-4+3i)^5}.
\end{equation}
Hypergeometric identities are:
\begin{eqnarray}
\!\hpg{2}{1}{\!3/20,\,7/20}{3/4}{\,x} &\equal&
{\textstyle\left(1-\frac{8(4+3i)}{25}x\right)}^{-1/8}
\hpg{2}{1}{1/40,\,9/40}{3/4}{\varphi_5(x)},\\
\hpg{2}{1}{3/20,\,2/5}{4/5}{\,x} &\equal & (1-x)^{-1/40}
{\textstyle\left(1-\frac{2+11i}4x\right)}^{-1/10}
\hpg{2}{1}{\!1/40,11/40}{4/5}{\frac1{\varphi_5(1/x)}}.\nonumber\\
\end{eqnarray}
\\

{\bf Acknowledgements.} The author would like to thank Robert S. Maier,
Frits Beukers and Masaaki Yoshida for useful references and remarks.


\end{document}